%% file: author.tex
\documentclass[graybox]{svmult}

\usepackage{mathptmx}       
\usepackage{helvet}         
\usepackage{courier}        
\usepackage{type1cm}        
\usepackage{makeidx}         
\usepackage{graphicx}        
\usepackage{multicol}        
\usepackage[bottom]{footmisc}

\makeindex             

\usepackage[font={footnotesize,it}]{caption}
\usepackage{amsmath}
\usepackage{hyperref}
\usepackage[noabbrev]{cleveref}
\usepackage[compact]{titlesec}
\usepackage{amssymb}
\usepackage{subfig}
\usepackage{pgfplots}
\pgfplotsset{compat=newest}

\begin{document}

\title*{A hybrid DEIM and leverage scores based method for CUR index selection}
\author{Perfect Y. Gidisu and Michiel E. Hochstenbach}
\institute{TU Eindhoven, PO Box 513, 5600MB, The Netherlands, \email{p.gidisu@tue.nl, m.e.hochstenbach@tue.nl}}

%
%

\maketitle

\vspace{-1mm}
\abstract*{The discrete empirical interpolation method (DEIM) may be used as an index selection strategy for formulating a CUR factorization. A notable drawback of the original DEIM algorithm is that the number of column or row indices that can be selected is limited to the number of input singular vectors.
We propose a new variant of DEIM, which we call L-DEIM, a combination of the strength of deterministic leverage scores and DEIM.  This method allows for the selection of a number of indices  greater than the number of available singular vectors. Since DEIM requires singular vectors as input matrices, L-DEIM is particularly attractive for example in big data problems when computing a rank-$k$-SVD approximation is expensive even for moderately small $k$ since it uses a lower-rank SVD approximation instead of the full rank-$k$ SVD. We empirically demonstrate the performance of L-DEIM, which despite its efficiency, may achieve comparable results to the original DEIM and even better approximations than some state-of-the-art methods.}

\abstract{The discrete empirical interpolation method (DEIM) may be used as an index selection strategy for formulating a CUR factorization. A notable drawback of the original DEIM algorithm is that the number of column or row indices that can be selected is limited to the number of input singular vectors.
We propose a new variant of DEIM, which we call L-DEIM, a combination of the strength of deterministic leverage scores and DEIM.  This method allows for the selection of a number of indices  greater than the number of input singular vectors. Since DEIM requires singular vectors as input matrices, L-DEIM is particularly attractive for example in big data problems when computing a rank-$k$ SVD approximation is expensive even for moderately small $k$ since it uses a lower-rank SVD approximation instead of the full rank-$k$ SVD. We empirically demonstrate the performance of L-DEIM, which despite its efficiency, may achieve comparable results to the original DEIM and even better approximations than some state-of-the-art methods.
}

\section{Introduction}
\label{sec:1}
Data sets are often represented by large matrices. In recent times, with the growth of the internet (industrial) data matrices are big and may be hard to manage. Examples of such data sets include text documents, customer databases, stocks, and financial transactions. In many data analyses, we need dimension reduction and for many applications, we need interpretable dimension reduction of which a CUR decomposition is one form. A CUR factorization is a low-rank matrix approximation proposed as an alternative to the TSVD to ensure interpretability and preserve relevant properties like sparsity or nonnegativity of the underlying matrix. A rank-$k$ CUR decomposition of an $m\times n$ matrix $A$ has the form
\vspace{-2mm}
\[
A \approx CMR = AP \cdot M \cdot S^ TA,
\]
where $C \in {\mathbb R}^{m\times k}$ and $R \in {\mathbb R}^{k\times n}$ are subsets of the columns and rows of $A$, respectively. The matrices $P\in {\mathbb R}^{n\times k}$ and $S\in {\mathbb R}^{m\times k}$ are index selection matrices with some columns of the identity indicating the columns and rows that are picked. The matrix $M$ is constructed to minimize the approximation error. There are several variants of this decomposition, which implies that the three factors are not necessarily unique. In \cite{Sorensen,Mahoney} the authors present algorithms for a CUR factorization based on a rank-$k$ singular value decomposition. Sorensen and Embree \cite{Sorensen} propose a CUR approximation using a discrete interpolation method (DEIM) on the rank-$k$ singular vectors.  The index selection method DEIM has first been introduced in the context of model order reduction \cite{Chaturantabut}.  In \cite{Sorensen}, it is shown to be a viable index selection method for identifying the most representative and influential subset of columns and rows that define a low-dimensional space of the data. The DEIM-induced CUR requires the computation of the SVD or its approximation. A notable limitation of this index selection algorithm is that the number of indices that can be selected is limited to the number of available singular vectors. In an attempt to address this, we propose a new extension called L-DEIM. The L-DEIM scheme combines the strengths of deterministic leverage scores sampling \cite{papailiopoulos2014} and the DEIM procedure. Our new approach is an alternative index selection method that is particularly attractive in a setting (for example big data problems) where we want a rank-$k$ CUR decomposition and computing a rank-$k$ SVD approximation is expensive even for moderately small $k$. This new algorithm allows us to select $k$ indices without having to compute the full rank-$k$ SVD by using a lower-rank SVD approximation instead. It may be viewed as an approach to reuse the same information to further improve the approximation.

We denote 2-norm by $\lVert{\cdot}\rVert$. We use MATLAB notation to index vectors and matrices; thus, $A(:,p)$ denotes the $k$ columns of $A$ whose corresponding indices are in vector $p \in {\mathbb N}_+^k$.

\section{Related Works}
\label{sec:2}
In this section, we briefly review some state-of-the-art deterministic algorithms for a CUR decomposition. These algorithms have been developed for the column subset selection problem or interpolative decomposition, but can be generalized for a CUR decomposition. We derive our proposed algorithm L-DEIM by combining two of the algorithms.

\subsection{Standard DEIM}\label{sec:21}
The DEIM is a discrete variant of the empirical interpolation method for approximating systems of nonlinear ordinary differential equations. In a recent paper by Sorensen and Embree \cite{Sorensen}, the authors use this method in the formulation of a CUR decomposition. The DEIM algorithm requires a full rank-$k$ SVD of $A$ to select at most $k$ column and or row indices of $A$.
 To illustrate how the indices are selected via the DEIM index selection method, we first define a projector which the authors in \cite{Chaturantabut} called an interpolatory projector. Suppose we want to preserve $k$ rows of $A$ and we have the rank-$k$ approximation of $A$ as
\begin{equation*}
\vspace{-1mm}
\begin{array}{ccccc}
    A&\approx & U && F , \\[-0.5mm]
    m\times n & & m\times k && k\times n
\end{array}
\vspace{-1mm}
\end{equation*}
where $U$ contains the top $k$ left singular vectors. The matrix $F$ is a coefficient matrix to be defined such that the above approximation preserves exactly the desired $k$ rows of $A$. Let ${\mathbf s}\in {\mathbb N}_+^k$ be an index vector with unique entries from the row index set \{$1,\dots,m$\} of $A$. Now let $S\in {\mathbb R}^{m\times k}$ be an index selection matrix  with some columns of the identity matrix that selects certain rows of $A$, i.e., $S=I(:,{\mathbf s})$. Assuming we want to keep desired rows in ${\mathbf s}$ in the approximation, viz., $S^TA\approx S^T(UF)$.
If $S^TU$ is nonsingular, the coefficient matrix $F$ can be determined uniquely; $F=(S^TU)^{-1}S^T\!A$. This implies $A\approx U(S^TU)^{-1}S^T\!A=\mathbb{S}A$.
The operator $\mathbb{S}$ is the DEIM interpolatory projector, an oblique projector. The name interpolatory comes from the fact that the projected matrix $\mathbb{S}A$ matches $A$ in the ${\mathbf s}$ entries. Note that we can obtain a similar projector using the right singular vectors. The DEIM algorithm processes the left singular vectors sequentially starting with the first dominant singular vector. Each step considers the next singular vector to obtain the next index. The selected indices are used to compute the interpolatory projector $\mathbb{S}$. The next index is selected by removing the direction of the interpolatory projection in the previous vectors from the subsequent one and finding the index of the entry with the largest magnitude in the residual vector (for more details see \cite{Sorensen}).

In \cite{Drmac}, Drmac and Gugercin proposed the Q-DEIM; a variant of DEIM which runs a column pivoted QR factorization on the transposes of the right and left singular vectors to select the column and row indices, respectively.

\subsection{Deterministic Leverage Score Sampling}\label{sec:22}
Part of the new extension borrows an idea from the leverage scores of a matrix $A$, which is defined below. We denote the $i$th row of $V_k$ by $[V_k]_{i,:}$.
\begin{definition}\label{def:1}
Given a matrix $A\in {\mathbb R}^{m\times n}$ with rank$(A)\ge k$, let $V_k$ contain its $k$ leading right singular vectors. The rank-$k$ leverage score of the $i$th column of $A$ is 
\[\ell_i=\lVert{[V_k]_{i,:}}\rVert^2, \quad i=1,\dots,n.\]
\end{definition}
 The deterministic leverage score sampling procedure selects columns of $A$ corresponding to the indices of the largest leverage scores for a given $k$.  This deterministic column selection method proposed by Jolliffe \cite{jolliffe1972discarding} is one of the first column subset selection algorithms. The leverage score sampling algorithm can extract at least $k$ column indices of $A$ and the upper bound on the number of indices that can be selected is not immediate. For more details on the algorithm and the bound on the number of columns to be sampled see, \cite[Sect.~3.1]{papailiopoulos2014}.

\section{L-DEIM}\label{subsec:2.1}  We now introduce the new extension of DEIM. Our starting point is the method from the earlier work \cite{Sorensen}, which derives a rank-$\widehat k$ CUR factorization by applying DEIM to the $\widehat k$ singular vectors. Given the promising results of this algorithm compared to other state-of-the-art methods for a CUR approximation, our proposed algorithm builds on the DEIM procedure. Constructing a rank-$\widehat k$ CUR decomposition using L-DEIM requires a rank-$k$ singular vectors where $\widehat k > k$. The integer $k$ is the number of available (approximate) singular vectors, while $\widehat k$ is the number of indices to be selected. To select the $\widehat k$ indices, the proposed method performs the original DEIM to find the first $k$ indices while keeping the residual singular vector in each index selection step of the DEIM procedure. The residual singular vector is the error between the input singular vector and its approximation from interpolating the previous singular vectors at the selected indices; as in line~2 of Algorithm~1\footnote{Note that the backslash operator used in the algorithm is a Matlab type notation for solving linear systems and least-squares problems.}. At the end of the iteration, using the idea of leverage scores, we compute the 2-norm of the rows of the residual singular vectors to select the additional $\widehat k - k$ indices. The procedure is summarized in Algorithm~1. Note that the vectors in $U$ in line~3 of Algorithm~1 are the residual singular vectors and not the original singular vectors.

\noindent\vrule height 0pt depth 0.5pt width \textwidth \\
{\bf Algorithm~1: L-DEIM index selection} \\[-2,5mm]
\vrule height 0pt depth 0.3pt width \textwidth \\
{\bf Input:} $U \in {\mathbb R}^{m \times k}$ and $V \in {\mathbb R}^{n \times k}$, target rank = $\widehat k$, with $k \le \widehat k\le \min(m,n)$\\
{\bf Output:} column and row indices ${\mathbf s}, {\mathbf p} \in {\mathbb N}_+^{ \widehat k}$, respectively, with non-repeating entries  \\
\begin{tabular}{ll}
& {\bf for} $j = 1, \dots, k$ \\
{\footnotesize 1:} & \phantom{M} ${\mathbf s}(j)$ =  $\text{argmax}_{1\le i\le m}~ |(U(:,j))_i|$ \\
{\footnotesize 2:} & \phantom{M} $U(:,j+1) = U(:,j+1)-U(:,~1:j)\cdot (U({\mathbf s},1:j)
\ \backslash \ U({\mathbf s},j+1))$ \\
{\footnotesize 3:} & Compute $\ell_i = \lVert{[U]_{i:}}\rVert$ \quad for $i=1,\dots, m$; sort $\ell$ in non-increasing order\\

{\footnotesize 4:} & Remove entries in  $\ell$ corresponding to the indices in ${\mathbf s}$\\
{\footnotesize 5:} & ${\mathbf s}' = \widehat k- k$ indices corresponding to $\widehat k-k$ largest entries of $\ell$ \\
{\footnotesize 6:} & ${\mathbf s} = [{\mathbf s};~{\mathbf s}']$\\
{\footnotesize 7:} & Perform {\footnotesize 1--6} on $V$ to get index set ${\mathbf p}$
\end{tabular} \\
\vrule height 0pt depth 0.5pt width \textwidth

From Algorithm~1, if $\widehat k = k$ then the algorithm reduces to the standard DEIM. We note that if the target rank is not specified, given $k$, we can select at least $k$ indices but the upper bound on the number of indices to be selected is not immediate; we can select an arbitrary number of indices. Similar to leverage scores sampling, the L-DEIM allows for oversampling of columns and or rows.

{\bf Error bounds.} 
Let us consider a fixed matrix $A\in {\mathbb R}^{m\times n}$ with rank $\rho \le \min(m,n)$. For an arbitrary $k$ with $1\le k\le \rho$, the best rank-$k$ approximation of $A$ $(A_k)$ provided by the SVD gives $\lVert{A-A_k}\rVert=\sigma_{k+1}(A)$ where $\sigma_{k+1}$ is the $(k+1)$st singular value of $A$. Suppose that we have a known target rank $k<\min(m,n)$, a good rank-$\widehat k$ approximation $A_{\widehat k}$ gives $\lVert{A-A_{\widehat k}}\rVert\le \tau \ \lVert{A-A_k}\rVert$, where $\tau >0$ is a modest tolerance and $k\le \widehat k \le r$ is the rank of the decomposition with oversampling. The following result unifies the theoretical bound results for $\lVert{A-CMR}\rVert$ in \cite[Sect.~4]{Sorensen} and \cite[Append.~1]{Hendryx2021}.

\begin{proposition}(See \cite[Sect.~4]{Sorensen}, \cite[Append.~1]{Hendryx2021})
Given $A\in {\mathbb R}^{m\times n}$ and $1\le k \le \widehat k\le \min(m,n)$, let $S\in {\mathbb R}^{m\times \widehat k}$, $P\in {\mathbb R}^{n\times \widehat k}$ be index selection matrices and the top $k$ left and $right$ singular vectors be $U\in {\mathbb R}^{m\times k}$ and $V\in {\mathbb R}^{n\times k}$, respectively. Let $C=AP \in {\mathbb R}^{m\times \widehat k}$ and $R=S^TA \in {\mathbb R}^{n\times \widehat k} $ be of full rank, assuming we compute $M$ as $(C^T C)^{-1}C^T A \ R^T\!(R R^T )^{-1}$ and $S^TU$ and $V^TP$ are of full rank we have 
\[\lVert{A-CMR}\rVert=[\sigma_{\min}^{-1}(V^TP) +\sigma_{\min}^{-1}(S^TU)]\cdot\sigma_{k+1}.\]
\end{proposition}
The above error bounds suggest the index selection method which minimizes the quantities $\sigma_{\min}^{-1}\ (V^TP)$ and $\sigma_{\min}^{-1}\ (S^TU)$ are theoretically desirable.

\section{Experiments}\label{sec:4}
We perform some experiments to compare the approximation quality and runtimes of the new method L-DEIM with the existing deterministic methods discussed in \cref{sec:2}. We use the relative error $\lVert{A-CMR}\rVert/\lVert{A}\rVert$ and runtimes for selecting the column and row indices as the evaluation criteria. Note that the runtimes reported here do not include the time for computing the singular vectors. We run the algorithms on three real data sets used in \cite{Sorensen, Mahoney}. The application domains of the data sets are Internet term document analysis, genetics, and collaborative filtering. The Internet term document data is from the Technion Repository of Text Categorization Datasets (TechTC). We use test 26, which consists of a collection of 139 documents on two topics with 15210 terms describing each document \cite{gabrilovich}. As in \cite{Sorensen}, the $139 \times 15210$ TechTC matrix rows are scaled to have a unit 2-norm. We take the cancer genetics data set GSE10072 from National Institutes of Health. This data set has 107 patients described by 22283 probes. There are 58 patients with tumors and 49 without. We center the $22283 \times 107$ genetics data matrix by subtracting the mean of each row from the entries in that row. The final data set is the Jester joke data set \cite{eigentaste}, which is often used as a benchmark for recommender system research. The data matrix consists of 73421 users and their ratings for 100 jokes. We only consider users who have ratings for all 100 jokes. We center the resulting $14116 \times 100$ matrix by subtracting the mean of each column from all entries in that column.

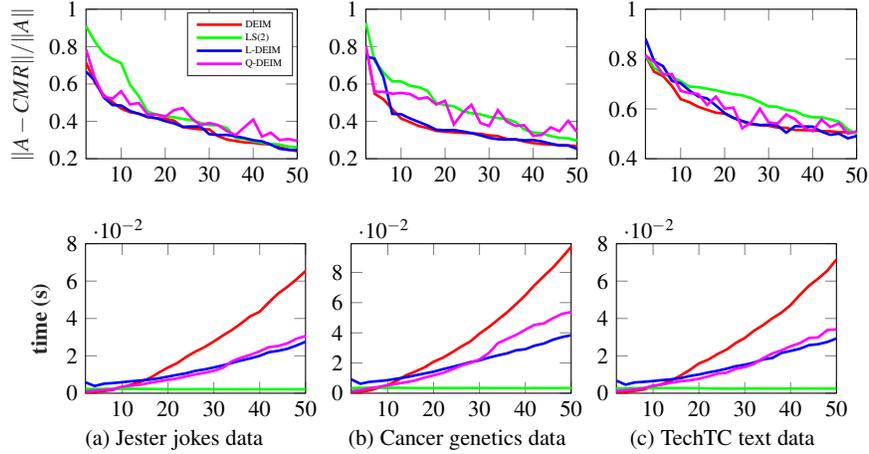
\begin{figure}[ht!]
\captionsetup[subfloat]{parskip=1pt,captionskip=1pt}
    \centering
    {\label{fig: a}{\input{jester_error}}}
	{\label{fig: b}{\input{Tumor_error}}}
	{\label{fig: c}{\input{Text_error}}}
	
    \subfloat[Jester jokes data]{\label{fig: 2a}{\input{jester_time}}}
	\subfloat[Cancer genetics data]{\label{fig: 2b}{\input{Tumor_time}}}
	\subfloat[TechTC text data]{\label{fig: 2c}{\input{Text_time}}}
	\caption{The approximation quality (first row) and runtimes (second row) of the L-DEIM scheme compared with the standard DEIM, Q-DEIM, and leverage scores sampling techniques using the three real data sets. Displayed are the relative errors $\lVert{A- CMR }\rVert/\lVert{A}\rVert$ and runtimes as a function of rank $k$.}\label{fig: 1}
\end{figure}

From \cref{fig: 1}, we see that the approximation quality of the proposed method L-DEIM can be as good as the original DEIM while the L-DEIM enjoys favorable runtimes. Both DEIM and L-DEIM have considerably lower approximation error than the other methods. The leverage scores sampling using two singular vectors seems to be the most efficient; however, we note that there is a trade-off between the runtimes and approximation quality. We show results of the leverage scores method using only the leading two singular vectors since higher choices yield worse approximation results.

\section{Conclusions}\label{sec:5}
We have presented a new extension of the DEIM index selection algorithm (L-DEIM) to identify additional indices for constructing a rank-$\widehat k$ CUR decomposition using a lower-rank SVD approximation. This is especially useful in a setting (for example big data problems) where computing a full rank-$\widehat k$ SVD is relatively expensive. The algorithm may be viewed not only as an extension of DEIM but also as an alternative index selection method for a CUR factorization. The L-DEIM procedure may also be suitable for point selection in the context of model order for nonlinear dynamical systems. Although the proposed algorithm is computationally more efficient than the original DEIM, experiments show that the approximation accuracy of both methods may be comparable when the target rank $\widehat k$ is at most twice the available $k$ singular vectors. For all results presented in \cref{sec:4}, we assume that given a target rank $\widehat k$, $2k=\widehat k$ in Algorithm~1. From experiments not presented here, if $\widehat k > 2k$ in Algorithm~1, then the rank-$\widehat k$ CUR approximation quality of the L-DEIM procedure which uses $k$ singular vectors may generally be worse than the rank-$\widehat k$ CUR factorization quality of the standard DEIM scheme which requires $\widehat k$ singular vectors. However, we stress that the L-DEIM is considerably cheaper.
A code for L-DEIM is available on \href{https://github.com/perfectyayra/L-DEIM-index-selection}{github.com/perfectyayra/L-DEIM-index-selection}.

{\bf Acknowledgement:}
This work has received funding from the European Union's Horizon 2020 research and innovation programme under the Marie Sk\l odowska-Curie grant agreement No 812912.

\input{referenc}
\end{document}

%% file: jester_error.tex
%
%
\definecolor{mycolor1}{rgb}{1.00000,0.00000,1.00000}%
\begin{tikzpicture}

\begin{axis}[%
width=0.24\columnwidth,
height=0.17\textwidth,
at={(0\textwidth,0\textwidth)},
scale only axis,
xmin=2,
xmax=50,
ymin=0.2,
ymax=1,
xtick={10,20,...,50},
ylabel style={font=\bfseries\color{white!15!black}},
ylabel={$\|A-C M R\|/ \|A\|$},
axis background/.style={fill=white},
legend style={legend columns=1, legend cell align=left, align=left, draw=white!15!black,nodes={scale=0.45, transform shape}}
]
\addplot [color=red, line width=1.0pt]
  table[row sep=crcr]{%
2	0.715930344143375\\
4	0.607794807540045\\
6	0.540767098731805\\
8	0.504038870363082\\
10	0.468121361706147\\
12	0.448972262712755\\
14	0.441684097388023\\
16	0.421531745894074\\
18	0.416793812204315\\
20	0.413765688283733\\
22	0.403897502848903\\
24	0.369410691135465\\
26	0.362304584214662\\
28	0.357514487764661\\
30	0.357194143256654\\
32	0.321074261599638\\
34	0.302858568350172\\
36	0.294832402321231\\
38	0.287253856915195\\
40	0.284994316185363\\
42	0.279153486813803\\
44	0.278177033914205\\
46	0.262108245032918\\
48	0.249426582168198\\
50	0.240610308027583\\
};
\addlegendentry{DEIM}

\addplot [color=green, line width=1.0pt]
  table[row sep=crcr]{%
2	0.910695316642416\\
4	0.827470495904533\\
6	0.763896345063271\\
8	0.734265071559396\\
10	0.711225170018758\\
12	0.596400930434316\\
14	0.534172913034389\\
16	0.451870294283042\\
18	0.436677263411425\\
20	0.42854421721496\\
22	0.422088568695811\\
24	0.41295793496364\\
26	0.406592990767729\\
28	0.392179182931688\\
30	0.379564984382633\\
32	0.378976213814182\\
34	0.364936380033347\\
36	0.310897274826454\\
38	0.306542530106533\\
40	0.293124639385017\\
42	0.285875397684023\\
44	0.279502321428433\\
46	0.275810074766779\\
48	0.265296826042405\\
50	0.26069110582709\\
};
\addlegendentry{LS(2)}

\addplot [color=blue, line width=1.0pt]
  table[row sep=crcr]{%
2	0.665997696726861\\
4	0.618700087215283\\
6	0.52882915259849\\
8	0.487462691588178\\
10	0.484737564183088\\
12	0.455302269709073\\
14	0.440210681636115\\
16	0.426024740135959\\
18	0.416901116975101\\
20	0.401120501220039\\
22	0.389424041558913\\
24	0.370764112704069\\
26	0.372582846088007\\
28	0.372537706717344\\
30	0.331593272094028\\
32	0.326858177368752\\
34	0.327455058327687\\
36	0.320433268025079\\
38	0.310733420719836\\
40	0.299112822890529\\
42	0.29319070922478\\
44	0.275277497415746\\
46	0.253416802155866\\
48	0.247156162540791\\
50	0.245687132796618\\
};
\addlegendentry{L-DEIM}

\addplot [color=mycolor1, line width=1.0pt]
  table[row sep=crcr]{%
2	0.786400758581953\\
4	0.633448388429203\\
6	0.532105832611197\\
8	0.521303665361497\\
10	0.561089949783087\\
12	0.490674186411081\\
14	0.498039790976647\\
16	0.437037434659745\\
18	0.436693783187423\\
20	0.424616594957831\\
22	0.46030747194791\\
24	0.471415293129027\\
26	0.417286425972753\\
28	0.37956910042685\\
30	0.391096985236472\\
32	0.380030441662481\\
34	0.348304106624483\\
36	0.326838256504742\\
38	0.36866581104114\\
40	0.409801611583529\\
42	0.318184187625146\\
44	0.339538840803639\\
46	0.299682297430404\\
48	0.305270957865858\\
50	0.295864175540893\\
};
\addlegendentry{Q-DEIM}

\end{axis}

\end{tikzpicture}%

%% file: Tumor_error.tex
%
%
\definecolor{mycolor1}{rgb}{1.00000,0.00000,1.00000}%
\begin{tikzpicture}

\begin{axis}[%
width=0.24\columnwidth,
height=0.17\textwidth,
at={(0\textwidth,0\textwidth)},
scale only axis,
xmin=2,
xmax=50,
ymin=0.2,
ymax=1,
xtick={10,20,...,50},
axis background/.style={fill=white}
]
\addplot [color=red, forget plot,line width=1.0pt]
  table[row sep=crcr]{%
2	0.803132278407496\\
4	0.547846120290618\\
6	0.517917681496192\\
8	0.465773515107938\\
10	0.415134893045559\\
12	0.396099818762933\\
14	0.375273960097903\\
16	0.364200418498391\\
18	0.3489520940572\\
20	0.344984709349083\\
22	0.341411168928289\\
24	0.338300998337735\\
26	0.336366693836066\\
28	0.333227474303734\\
30	0.324744003936846\\
32	0.321038257413838\\
34	0.305691634958101\\
36	0.304459999839793\\
38	0.291984182409046\\
40	0.283674939699599\\
42	0.278616546851746\\
44	0.275400559661381\\
46	0.273728703376462\\
48	0.270964297941545\\
50	0.268132135895247\\
};
\addplot [color=green, forget plot,line width=1.0pt]
  table[row sep=crcr]{%
2	0.927156970806745\\
4	0.7227118149332\\
6	0.661137656229397\\
8	0.614020695445062\\
10	0.613823992452984\\
12	0.591230229520993\\
14	0.583084273692498\\
16	0.563727538839352\\
18	0.49856041049882\\
20	0.486880436853214\\
22	0.481048200081526\\
24	0.454126270797648\\
26	0.443929920405985\\
28	0.43735644386481\\
30	0.424736605405632\\
32	0.417464924933014\\
34	0.412029223533785\\
36	0.396734215337619\\
38	0.356355355555573\\
40	0.341030262630174\\
42	0.336413996510631\\
44	0.321128468385939\\
46	0.312266528210305\\
48	0.309404990779138\\
50	0.296161660226848\\
};
\addplot [color=blue, forget plot,line width=1.0pt]
  table[row sep=crcr]{%
2	0.746937513281335\\
4	0.735849897425823\\
6	0.63553137519131\\
8	0.440158003842208\\
10	0.438359338492594\\
12	0.415999097601827\\
14	0.396955868652566\\
16	0.378345872075763\\
18	0.355255231602981\\
20	0.352442654733314\\
22	0.353450364035935\\
24	0.346205795051795\\
26	0.337394777225081\\
28	0.325126064304024\\
30	0.31181910626607\\
32	0.304620933950276\\
34	0.302512617344616\\
36	0.302798027670946\\
38	0.298398934650143\\
40	0.299990184073794\\
42	0.292548483789561\\
44	0.280800527442346\\
46	0.270902321656305\\
48	0.27082898753517\\
50	0.252452486733873\\
};
\addplot [color=mycolor1, forget plot,line width=1.0pt]
  table[row sep=crcr]{%
2	0.803132278407496\\
4	0.558546859027752\\
6	0.558373931187649\\
8	0.547712972919062\\
10	0.551907908222708\\
12	0.546548775795619\\
14	0.520655409966042\\
16	0.528604864861882\\
18	0.490508021442312\\
20	0.510133125129258\\
22	0.382943946398883\\
24	0.446103517683429\\
26	0.488545200437812\\
28	0.392031652149972\\
30	0.375141570265414\\
32	0.460537319381415\\
34	0.399273061954008\\
36	0.375799596121105\\
38	0.379275456448785\\
40	0.322732105615987\\
42	0.326997863663654\\
44	0.367740156541169\\
46	0.34725083965897\\
48	0.403080234392757\\
50	0.34222337395847\\
};

\end{axis}

\end{tikzpicture}%

%% file: Text_error.tex
%
%
\definecolor{mycolor1}{rgb}{1.00000,0.00000,1.00000}%
\begin{tikzpicture}

\begin{axis}[%
width=0.24\columnwidth,
height=0.17\textwidth,
at={(0\textwidth,0\textwidth)},
scale only axis,
xmin=2,
xmax=50,
ymin=0.4,
ymax=1,
xtick={10,20,...,50},
axis background/.style={fill=white}
]
\addplot [color=red, forget plot,line width=1.0pt]
  table[row sep=crcr]{%
2	0.814973528835639\\
4	0.748635797430956\\
6	0.731944734163647\\
8	0.692398487456192\\
10	0.639522160883684\\
12	0.626569202998987\\
14	0.608338727510651\\
16	0.596338829116841\\
18	0.584868036170298\\
20	0.580782901581641\\
22	0.563628843302438\\
24	0.549952166946957\\
26	0.537856739918729\\
28	0.535245902803659\\
30	0.532279712660795\\
32	0.525097954456463\\
34	0.523804275272877\\
36	0.515475835837211\\
38	0.514602648160535\\
40	0.513508127546361\\
42	0.511053413080444\\
44	0.510129353898988\\
46	0.506527911456114\\
48	0.505994584004329\\
50	0.505725485174587\\
};
\addplot [color=green, forget plot,line width=1.0pt]
  table[row sep=crcr]{%
2	0.814973528835639\\
4	0.764167972954265\\
6	0.738936978169535\\
8	0.728717552154486\\
10	0.711010961661154\\
12	0.689591013509474\\
14	0.687222594330708\\
16	0.681141756996214\\
18	0.674520421417697\\
20	0.66624789834417\\
22	0.659154564785967\\
24	0.654160438237681\\
26	0.643780127527925\\
28	0.627343414696671\\
30	0.610081402592877\\
32	0.608761053604095\\
34	0.595495094394222\\
36	0.591868412655947\\
38	0.57415497670665\\
40	0.566395028273817\\
42	0.566034710994871\\
44	0.558277666783151\\
46	0.548804379564007\\
48	0.518682200185233\\
50	0.498625648567075\\
};
\addplot [color=blue, forget plot,line width=1.0pt]
  table[row sep=crcr]{%
2	0.882654596166051\\
4	0.790496146572939\\
6	0.771677297892373\\
8	0.712826766946618\\
10	0.702667174581533\\
12	0.673480255646292\\
14	0.650191511201525\\
16	0.641964907683613\\
18	0.613639307813661\\
20	0.583292800431394\\
22	0.564080657517228\\
24	0.551452397602077\\
26	0.541016740213403\\
28	0.535064215650531\\
30	0.536528104377513\\
32	0.533026032108363\\
34	0.504704557111881\\
36	0.530636706831853\\
38	0.529488007572838\\
40	0.528717944759116\\
42	0.508931481858104\\
44	0.496036646146593\\
46	0.500629570312379\\
48	0.481752883569875\\
50	0.491408365199589\\
};
\addplot [color=mycolor1, forget plot,line width=1.0pt]
  table[row sep=crcr]{%
2	0.814973528835639\\
4	0.7876871693996\\
6	0.739030703667287\\
8	0.739842073073218\\
10	0.672793219844087\\
12	0.663268573376545\\
14	0.659826670020707\\
16	0.615974428656283\\
18	0.64943144972648\\
20	0.600830170543645\\
22	0.605106778316872\\
24	0.521792412202962\\
26	0.545752643580124\\
28	0.600569298973642\\
30	0.546136941048532\\
32	0.540784310249314\\
34	0.577531684298316\\
36	0.558402221852793\\
38	0.525574437641155\\
40	0.510256870873859\\
42	0.557822059552413\\
44	0.527177978636054\\
46	0.541512471273237\\
48	0.500634245223173\\
50	0.510247517856023\\
};

\end{axis}

\end{tikzpicture}%

%% file: jester_time.tex
%
%
\definecolor{mycolor1}{rgb}{1.00000,0.00000,1.00000}%
\begin{tikzpicture}

\begin{axis}[%
width=0.25\columnwidth,
height=0.17\textwidth,
at={(0cm,0cm)},
scale only axis,
xmin=2,
xmax=50,
ymin=0,
ymax=0.08,
xtick={10,20,...,50},
ylabel style={font=\bfseries\color{white!15!black}},
ylabel={time (s)},
axis background/.style={fill=white}
]
\addplot [color=red, forget plot, line width=1.0pt]
  table[row sep=crcr]{%
2	0.000261738\\
4	0.000629641\\
6	0.001160668\\
8	0.001842615\\
10	0.003026179\\
12	0.004624341\\
14	0.006108643\\
16	0.00806242\\
18	0.011033919\\
20	0.013762086\\
22	0.016087955\\
24	0.019473361\\
26	0.022307243\\
28	0.02473052\\
30	0.02772309\\
32	0.030924824\\
34	0.033927256\\
36	0.037463429\\
38	0.041279926\\
40	0.043612136\\
42	0.048644262\\
44	0.053301916\\
46	0.056857583\\
48	0.060920498\\
50	0.065379053\\
};
\addplot [color=green, forget plot, line width=1.0pt]
  table[row sep=crcr]{%
2	0.002266837\\
4	0.002225521\\
6	0.002279109\\
8	0.002275226\\
10	0.002328129\\
12	0.002354083\\
14	0.00232434\\
16	0.002293166\\
18	0.002220175\\
20	0.002147851\\
22	0.002147227\\
24	0.002145467\\
26	0.002147696\\
28	0.002189053\\
30	0.002153834\\
32	0.002133891\\
34	0.00213583\\
36	0.002114678\\
38	0.002139794\\
40	0.002139491\\
42	0.002113479\\
44	0.002117368\\
46	0.002129781\\
48	0.002094243\\
50	0.002131303\\
};
\addplot [color=blue, forget plot, line width=1.0pt]
  table[row sep=crcr]{%
2	0.005782413\\
4	0.003864077\\
6	0.005091795\\
8	0.005440194\\
10	0.00589494\\
12	0.006400627\\
14	0.006798365\\
16	0.007265928\\
18	0.007970007\\
20	0.008894443\\
22	0.009832831\\
24	0.010688721\\
26	0.011980378\\
28	0.012749606\\
30	0.013802206\\
32	0.014910569\\
34	0.016244827\\
36	0.017547031\\
38	0.018682023\\
40	0.020017375\\
42	0.021951161\\
44	0.022672895\\
46	0.023757563\\
48	0.025522322\\
50	0.027570026\\
};
\addplot [color=mycolor1, forget plot, line width=1.0pt]
  table[row sep=crcr]{%
2	0.00089383\\
4	0.001278162\\
6	0.001647346\\
8	0.00224147\\
10	0.003515088\\
12	0.004074082\\
14	0.004569895\\
16	0.005533955\\
18	0.00637847\\
20	0.007152052\\
22	0.007885468\\
24	0.009043245\\
26	0.010068482\\
28	0.010935945\\
30	0.0119717\\
32	0.013664765\\
34	0.017162838\\
36	0.019070083\\
38	0.020522713\\
40	0.022451005\\
42	0.024663283\\
44	0.025281075\\
46	0.026628105\\
48	0.029146173\\
50	0.030606961\\
};
\end{axis}
\end{tikzpicture}%

%% file: Tumor_time.tex
%
%
\definecolor{mycolor1}{rgb}{1.00000,0.00000,1.00000}%
\begin{tikzpicture}

\begin{axis}[%
width=0.25\columnwidth,
height=0.17\textwidth,
at={(0cm,0cm)},
scale only axis,
xmin=2,
xmax=50,
ymin=0,
ymax=0.099,
xtick={10,20,...,50},
axis background/.style={fill=white}
]
\addplot [color=red, forget plot,line width=1.0pt]
  table[row sep=crcr]{%
2	0.000303037\\
4	0.000875399\\
6	0.001850977\\
8	0.00346837\\
10	0.005384108\\
12	0.008315872\\
14	0.011889297\\
16	0.01438239\\
18	0.017112779\\
20	0.020873998\\
22	0.023519671\\
24	0.027011683\\
26	0.030882763\\
28	0.0342839\\
30	0.039441095\\
32	0.043752218\\
34	0.048619848\\
36	0.053513096\\
38	0.059039237\\
40	0.064905355\\
42	0.071525935\\
44	0.077450848\\
46	0.083177894\\
48	0.089870842\\
50	0.096987178\\
};
\addplot [color=green, forget plot,line width=1.0pt]
  table[row sep=crcr]{%
2	0.003365769\\
4	0.003553765\\
6	0.003512787\\
8	0.00364937\\
10	0.003543518\\
12	0.003642748\\
14	0.003482517\\
16	0.003421649\\
18	0.003369365\\
20	0.003379546\\
22	0.003445812\\
24	0.00341251\\
26	0.00341815\\
28	0.003438537\\
30	0.003420777\\
32	0.003432273\\
34	0.003343001\\
36	0.003353583\\
38	0.003385541\\
40	0.003328563\\
42	0.003456355\\
44	0.003348112\\
46	0.003387819\\
48	0.003396628\\
50	0.003347929\\
};
\addplot [color=blue, forget plot,line width=1.0pt]
  table[row sep=crcr]{%
2	0.009196436\\
4	0.006246395\\
6	0.007501965\\
8	0.008029278\\
10	0.008633176\\
12	0.009647991\\
14	0.010569647\\
16	0.011823667\\
18	0.013062412\\
20	0.014268152\\
22	0.015894971\\
24	0.017388993\\
26	0.018904981\\
28	0.020569093\\
30	0.021739884\\
32	0.023560665\\
34	0.025051623\\
36	0.026389478\\
38	0.028444075\\
40	0.029151153\\
42	0.031495317\\
44	0.032849582\\
46	0.035415567\\
48	0.037164496\\
50	0.038394927\\
};
\addplot [color=mycolor1, forget plot,line width=1.0pt]
  table[row sep=crcr]{%
2	0.001320207\\
4	0.001899994\\
6	0.003296179\\
8	0.004146534\\
10	0.005282642\\
12	0.006276564\\
14	0.007281716\\
16	0.008476953\\
18	0.010329357\\
20	0.011876164\\
22	0.013384111\\
24	0.015698295\\
26	0.018296949\\
28	0.020111569\\
30	0.022378888\\
32	0.026510834\\
34	0.033119072\\
36	0.035941225\\
38	0.03864317\\
40	0.042094347\\
42	0.045237884\\
44	0.046154824\\
46	0.049619881\\
48	0.052517396\\
50	0.053791285\\
};
\end{axis}
\end{tikzpicture}%

%% file: Text_time.tex
%
%
\definecolor{mycolor1}{rgb}{1.00000,0.00000,1.00000}%
\begin{tikzpicture}

\begin{axis}[%
width=0.25\columnwidth,
height=0.17\textwidth,
at={(0cm,0cm)},
scale only axis,
xmin=2,
xmax=50,
ymin=0,
ymax=0.08,
xtick={10,20,...,50},
axis background/.style={fill=white}
]
\addplot [color=red, forget plot,line width=1.0pt]
  table[row sep=crcr]{%
2	0.000246703\\
4	0.000678981\\
6	0.001300719\\
8	0.002018298\\
10	0.003550519\\
12	0.004946557\\
14	0.006855544\\
16	0.00946175\\
18	0.012443046\\
20	0.015925987\\
22	0.017912309\\
24	0.020881149\\
26	0.023307473\\
28	0.026806418\\
30	0.029550122\\
32	0.03338906\\
34	0.036448283\\
36	0.040045239\\
38	0.043119569\\
40	0.04719361\\
42	0.052549669\\
44	0.057480389\\
46	0.061287268\\
48	0.065531633\\
50	0.071591327\\
};
\addplot [color=green, forget plot,line width=1.0pt]
  table[row sep=crcr]{%
2	0.002478138\\
4	0.002677543\\
6	0.00266977\\
8	0.00276626\\
10	0.002664163\\
12	0.002710331\\
14	0.002760659\\
16	0.002719973\\
18	0.002703434\\
20	0.00263588\\
22	0.002529405\\
24	0.002489654\\
26	0.002456915\\
28	0.002433454\\
30	0.00250611\\
32	0.002503088\\
34	0.002455075\\
36	0.002493597\\
38	0.002481763\\
40	0.002475163\\
42	0.002477621\\
44	0.002498486\\
46	0.002496167\\
48	0.002468684\\
50	0.002484027\\
};
\addplot [color=blue, forget plot,line width=1.0pt]
  table[row sep=crcr]{%
2	0.006698971\\
4	0.004582651\\
6	0.00575448\\
8	0.006088143\\
10	0.006523576\\
12	0.007019213\\
14	0.007440244\\
16	0.008063184\\
18	0.009213254\\
20	0.009994074\\
22	0.011105145\\
24	0.011988051\\
26	0.013083884\\
28	0.014443946\\
30	0.015251577\\
32	0.016349543\\
34	0.018095844\\
36	0.018834282\\
38	0.021564189\\
40	0.022635441\\
42	0.023777841\\
44	0.025702343\\
46	0.026224799\\
48	0.027398959\\
50	0.029279544\\
};
\addplot [color=mycolor1, forget plot,line width=1.0pt]
  table[row sep=crcr]{%
2	0.001124121\\
4	0.001571881\\
6	0.001983498\\
8	0.002643066\\
10	0.004060603\\
12	0.004602451\\
14	0.005715664\\
16	0.006394809\\
18	0.006928471\\
20	0.007734403\\
22	0.008988577\\
24	0.010431881\\
26	0.010967662\\
28	0.013016193\\
30	0.013985935\\
32	0.015696264\\
34	0.018718726\\
36	0.021226744\\
38	0.022906603\\
40	0.025102224\\
42	0.026775096\\
44	0.029550419\\
46	0.029865065\\
48	0.033764068\\
50	0.034182859\\
};
\end{axis}
\end{tikzpicture}%